\documentclass{nyj}
\title[Correction to Knotted Hamiltonian cycles]{Correction to ``Knotted Hamiltonian cycles in spatial embeddings of complete graphs"} % If the title is too long for the running heads, give a
\author [J.Foisy]{Joel Foisy}

\address {Department of Mathematics, SUNY Potsdam, Potsdam, NY 13676}

\email{foisyjs@potsdam.edu}



\usepackage{latexsym}
\usepackage{amsmath, amsthm}
\usepackage{graphicx}
\usepackage{amssymb}

\newtheorem{thm}{Theorem}[section]

\newtheorem{lem}[thm]{Lemma}

\begin{document} 
       % Abstracts are required, and will be distributed by a listserv list.
       % Please place the text of your abstract in the next environment.

\begin{abstract}
We state and prove a correct version of a theorem presented in \cite{Foisy}.
\end{abstract}
\maketitle
%\tableofcontents

%%%% **** The text of the paper starts here **** %%%%

Professor Masakazu Teragaito has recently pointed out that Theorem 3.3 of \cite{Foisy} is incorrect as stated.  The fact that $\mu_{f}(G,\Gamma;6)=3$ is 
independent of embedding of $K_8$ does not necessarily imply that there are at least $3$ knotted Hamiltonian cycles.  For example, there could be exactly one knotted Hamiltonian cycle with $a_2(K)=3$.

Professor Kouki Taniyama has further pointed out that Lemma 2, in \cite{shim}, has a gap in its proof.  He has proposed a rigorous proof of a weaker version Lemma 2.  In this short paper, we will state and prove this weaker version of Lemma 2 in \cite{shim}, and then apply it to obtain a weaker version of Theorem 3.3 of \cite{Foisy}.  For definitions of terms, see \cite{shim} and \cite{Foisy}.  Here is the modified version of Shimabara's Lemma 2 that we will prove:

\begin{lem} \textit{Let $\Gamma$ be a set of cycles in an undirected 
graph G.  The invariant $\mu$$_{f}(G,\Gamma;n)$ does not depend on 
the spatial embedding $f$ of G if the following two conditions hold}:
\label{shimlem}

(1)  For any edges A, B, E such that A is adjacent to B, \newline
$\nu_{1}(\Gamma;A,B,E)\equiv 0$     $(mod$ $n)$.

(2)  For any pairs of non-adjacent edges (A,B) and (E,F), \newline
 $\nu_{2}(\Gamma;A,B;E,F)\equiv 0$     $(mod$ $2n)$.
\label{new}
\end{lem}

The difference between this new lemma and the original comes in the second condition, where the equivalence is mod $2n$, not $n$.  For the proof of Lemma 2, case 2, on p. 410 of \cite{shim}, the definition of linking number used does not work.  For the version of linking number used, $\zeta(A,B)$ depends on the order of $A$ and $B$, whereas the equality $\displaystyle\sum_{E,F}\displaystyle\sum_{\gamma \in \Gamma_1} \epsilon(c)\zeta(f_{\gamma}(E),f_{\gamma}(F))=\displaystyle\sum_{E,F}(n_3-n_4) \zeta(f(E),f(F))$ implicitly uses the assumption that $\zeta(A,B)=\zeta(B,A)$. %put mod n in?

A proof of Lemma \ref{new} is possible if one uses a different (but equivalent) version of $\zeta$ and linking number.  For $A$ and $B$ two disjoint oriented arcs or circles in ${\mathbb R}^3$, define $\zeta'(A,B)=(1/2)\displaystyle\sum_{c}\epsilon (c)$, with the summation being over all crossings between $A$ and $B$.  If $A$ and $B$ are circles, then $\zeta'(A,B)$ gives the linking number of $A$ and $B$, $lk(A,B)$.  It then follows immediately that $\zeta'(A,B)=\zeta'(B,A)$.  The proof of Lemma \ref{new} proceeds as in the proof of Lemma 2 in \cite{shim}, but now ends with:  
$\delta(u)=\displaystyle\sum_{E,F}\displaystyle\sum_{\gamma \in \Gamma_1} \epsilon(c)\zeta(f_{\gamma}(E),f_{\gamma}(F))=\displaystyle\sum_{E,F}(n_3-n_4) \zeta'(f(E),f(F))$.%mod n

\noindent
It then follows, since each $\zeta'(f(E),f(F))$ is either an integer or an integer divided by $2$, that, $\delta(u)$ is congruent to $0$ mod $n$ if $|n_3-n_4| \equiv 0$ mod $2n$.

Now, in \cite{Foisy}, it was shown that $\nu_2 \equiv 0$ mod $6$.  It was also shown that there exists an embedding of $K_8$ with exactly $21$ knotted Hamiltonian cycles, each with Arf invariant $1$.  One can also verify that each of these knotted cycles is a trefoil, with $a_2=1$.  This embedding, together with Lemma \ref{new}, implies that $\mu_{f}(K_8,\Gamma, 3) \equiv 0$ for every spatial embedding of $K_8$.  By Theorem 2.2 of \cite {Foisy}, there is at least one Hamiltonian cycle with Arf invariant $1$, in every spatial embedding of $K_8$.

  We thus have the following corrected version of Theorem 3.3 from \cite {Foisy}.

\begin{thm}

Given an embedding of $K_8$, at least one of the following must occur in that embedding:

\begin{enumerate}

\item At least $3$ knotted Hamiltonian cycles.

\item  Exactly $2$ knotted Hamiltonian cycles $C_1$ and $C_2$, with $a_2(C_1) \equiv 1$ $($mod $3)$ and $a_2(C_2) \equiv 2$ $($mod $3)$, or $a_2(C_1) \equiv 0$ $($mod $3)$ and $a_2(C_2) \equiv 0$ $($mod $3)$.  Either $C_1$ or $C_2$ has non-zero Arf invariant.

\item  Exactly $1$ knotted Hamiltonian cycle, $C$ with $a_2(C) \equiv 1$ $($mod $2)$ and $a_2(C) \equiv 0$ $($mod $3)$.  (Equivalently:  $a_2(C) \equiv 3$ $($mod $6)$.)

\end{enumerate}

\end{thm}

It thus remains an open question to determine if $1$ is the best lower bound for the minimum number of knotted Hamiltonian cycles in every spatial embedding of $K_8$.

\bigskip
\noindent
{\bf Acknowlegments}  The author would like to thank Professors Masakazu Teragaito, Kouki Taniyama and Ryo Nikkuni for valuable comments and suggestions.

\end{document}